\documentclass[12pt]{amsart}
\usepackage{amsmath,amssymb,amsfonts,amscd,verbatim,graphicx}
\usepackage{xypic}
\newtheorem{thm}{Theorem}

\theoremstyle{definition}
\newtheorem{definition}[subsection]{Definition}

\theoremstyle{remark}

\numberwithin{equation}{section}
\newtheorem{cor}{Corollary}[section]

\theoremstyle{definition}
\newtheorem{defi}[cor]{Definition}

\newcommand{\co}{\colon\thinspace}

\begin{document}

\title[Robust $4-$manifolds]{Robust ${\mathbf 4-}$manifolds and robust embeddings}
\author{Vyacheslav Krushkal}
\address{Department of Mathematics, University of Virginia, Charlottesville, VA 22904}
\email{krushkal\char 64 virginia.edu}

\thanks{This research was supported in part by the NSF}

\begin{abstract}
A link in the $3-$sphere is homotopically trivial, according to Milnor,
if its components bound disjoint maps of disks in the $4-$ball.
This paper concerns the question of what spaces
give rise to the same class of homotopically trivial links when used in place of disks in an analogous definition.
We show that there are $4-$manifolds for which this property depends on their embedding in the $4-$ball.
This work is motivated by the A-B slice problem, a reformulation of the
$4-$dimensional topological surgery conjecture.
As a corollary this provides a new, secondary, obstruction in the A-B
slice problem for a certain class of decompositions of $D^4$.
\end{abstract}

\maketitle

\section{Introduction}
Classification of knots and links up to concordance, in both smooth and topological categories,
is an important and difficult problem in $4-$dimensional topology. Recall that
two links in the $3-$sphere are concordant if they bound disjoint embeddings
of (smooth or locally flat, depending on the category) annuli in $S^3\times [0,1]$.
In \cite{M} Milnor introduced the notion of {\em link homotopy}, often
referred to as the ``theory of links modulo knots'', which turned out to be
much more tractable. In particular, there is an elegant
characterization of homotopically trivial links using the {\em Milnor group},
a certain rather natural nilpotent quotient of the fundamental group of the link. Two links are link homotopic
if they bound disjoint maps of annuli in $S^3\times [0,1]$, so the annuli are disjoint
from each other but unlike the definition of concordance, they are allowed to have
self-intersections.

The subject of link homotopy brings together $4-$dimensional geometric topology and
the classical techniques of nilpotent group theory. An area where both approaches are important, and which is
a motivation for the results in this paper, is the {\em A-B slice problem}, a reformulation of the
topological $4-$dimensional surgery conjecture \cite{FL}, \cite{FQ}. Roughly speaking, this paper concerns the problem of characterizing
spaces (which in interesting cases are $4-$manifolds
with a specified curve in the boundary)  which give rise to the same class of homotopically trivial
links when used in place of disks in a definition analogous to Milnor's.

A specific question concerning the
A-B slice problem and related to link homotopy theory is the following. Suppose $M$ is
a codimension zero submanifold of the $4-$ball, $i\co M\hookrightarrow D^4$, with a specified curve
${\gamma}\subset \partial M$ forming a knot in the $3-$sphere: $i({\gamma})\subset S^3=\partial D^4$.
Given such a pair $(M,{\gamma})$, does there exist a homotopically essential link $(i_1({\gamma}),\ldots,i_n({\gamma}))$
in the $3-$sphere formed by disjoint embeddings $i_1,\ldots,i_n$ of $(M,{\gamma})$ into $(D^4,S^3)$? If the answer is negative, the
pair $(M,{\gamma})$ is called {\em robust}. The analysis of this problem is substantially
more involved than the classical link homotopy case where one considers disks with self-intersections:
in general the $1-$ and $2-$handles of $4-$manifolds embedded in
$D^4$ may link, and the relations in the fundamental group of the complement do not
have the ``standard'' form implied by the Clifford tori in Milnor's theory.

The main result of this paper is the existence of $4-$manifolds for which this property
depends on their embedding in the $4-$ball:

\begin{thm} \label{thm} \sl There exist submanifolds $i\co (M,{\gamma})\hookrightarrow (D^4,S^3)$
such that

(i) There are disjoint embeddings of several copies $(M_i, {\gamma}_i)$ of $(M,{\gamma})$ into $(D^4,S^3)$
forming a homotopically essential link $({\gamma}_1, \ldots, {\gamma}_n)$ in the $3-$sphere.

(ii) Given any disjoint embeddings of $(M,{\gamma})$ into $(D^4,S^3)$, each one
isotopic to the original embedding $i$, the link $({\gamma}_1, \ldots, {\gamma}_n)$ formed
by their attaching curves in the $3-$sphere is homotopically trivial.
\end{thm}

In the A-B slice problem one
considers {\em decompositions} of the $4-$ball, $D^4=A\cup B$, where the specified curves
${\alpha}, {\beta}$ of the two parts form the Hopf link in $S^3=\partial D^4$ (see \cite{FL}, \cite{K1}, and section
\ref{sec:robust} below). It is shown
in \cite{K1} that there exist decompositions where neither of the two sides is robust. This result left open the question of whether these decompositions, in fact, may be used  to solve the general $4-$dimensional topological surgery conjecture.

This paper provides a detailed analysis of the construction in \cite{K1}. The proof of theorem \ref{thm} shows
that in this case,
one of the two parts of the decomposition is robust provided that the re-embeddings forming the
link $({\gamma}_1,\ldots, {\gamma}_n)$ are topologically equivalent to the original embedding
into the $4-$ball.
This provides a new obstruction in the A-B slice problem for this class of
decompositions of the $4-$ball, exhibiting a new phenomenon
where the obstruction depends not just on the submanifold but also on its specific
embedding into $D^4$. In particular, this shows that the construction in \cite{K1} does not
satisfy the equivariance condition which is necessary for solving the canonical $4-$dimensional
surgery problems.
An important open question is whether given
{\em any} decomposition $D^4=A\cup B$, the conclusion (ii) of the theorem holds for either $A$ or $B$.

The main tool in the proof of part (ii) of the theorem is the Milnor group in the context of the {\em relative-slice problem}.
The context here
is substantially different from Milnor's original work, since we use it to analyze
embeddings of more general submanifolds in the $4-$ball, and the topology of these spaces is richer
than the setting of disks with self-intersections in the $4$-ball, considered in classical link homotopy. 
To find an obstruction, one has to consider in detail the structure of the graded Lie algebra 
associated to the lower central series of the link group.
The strategy of the proof
should be useful in further study of the A-B slice problem. 

Section \ref{sec:robust} gives a detailed definition of robust $4-$manifolds and robust embeddings,
and discusses its relation with the A-B slice problem. The construction \cite{K1} of the submanifolds,
used in the proof of theorem \ref{thm}, is recalled in section \ref{decompositions}. We review
the Milnor group in the $4-$dimensional setting and complete the
proof of theorem \ref{thm} in section \ref{sec:Milnor}.

\section{Robust 4-manifolds and the A-B slice problem} \label{sec:robust}

This section states the definition of robust $4-$manifolds and robust embeddings,
the notions which provide a convenient setting for the results
in this paper, and important in relation with the A-B slice problem.
Let $M$ be a $4-$manifold with a specified curve in its boundary, ${\gamma}\subset \partial M$.
Let $i\co M \hookrightarrow D^4$ be an embedding into the $4-$ball with $i({\gamma})\subset
S^3=\partial D^4$.

\begin{definition} The pair $(M,{\gamma})$ is {\em robust} if given any $n\geq 2$ and disjoint embeddings
$i_1, \ldots, i_n$ of $(M,{\gamma})$ into $(D^4,S^3)$, the link formed by the curves
$i_1({\gamma}),\ldots, i_n({\gamma})$ in the $3-$sphere
is homotopically trivial.

An {\em embedding} $i\co (M,{\gamma}) \hookrightarrow (D^4, S^3)$ is {\em robust} if
given any $n\geq 2$ and disjoint embeddings
$i_1, \ldots, i_n$ of $(M,{\gamma})$ into $(D^4,S^3)$, each one
isotopic to the original embedding $i$, the link formed by the curves
$i_1({\gamma}),\ldots, i_n({\gamma})$ in the $3-$sphere
is homotopically trivial. (In this case, we say that the re-embeddings are {\em standard}.)
\end{definition}

In these terms, theorem \ref{thm} states that there exist $(M,{\gamma})$ which are not robust
but which admit robust embeddings, first examples of this phenomenon.

It follows immediately from definition that the $2-$handle
$(D^2\times D^2, \{0\} \times \partial D^2)$, and more generally any kinky handle (a regular
neighborhood in the $4-$ball of a disk with self-intersections) is robust.

It is not difficult to give further examples: it follows from the link composition lemma
\cite{FL, KT} that the $4-$manifold $(B_o,{\beta})$ in figure \ref{fig:Bing pair}, obtained from the collar
${\beta}\times D^2\times [0,1]$ by attaching $2-$handles to the Bing double of the core of the
solid torus, is robust. This example illustrates the important point that the disjoint copies $i_j(M)$
in the definition above are {\em embedded}: it is easy to see that if the $2-$handles $H_1, H_2$ in figure
\ref{fig:Bing pair} were allowed to intersect, this $4-$manifold may be mapped to the collar on its attaching
curve and therefore there exist disjoint {\em singular} maps of copies of this manifold such that their
attaching curves $\{ {\gamma}_i \}$ form a homotopically essential
link in the $3-$sphere.

The complement in $D^4$ of the standard embedding of the $4-$manifold in figure \ref{fig:Bing pair}
is the $4-$manifold $A=$(genus one surface with one boundary component $\alpha$)$\times D^2$.
It is easy to see that $(A,{\alpha})$ is not robust: for example, the Borromean rings form a homotopically
essential link bounding disjoint standard genus one surfaces
in the $4-$ball.

To review the relation of these results with the $4-$dimensional topological surgery conjecture,
recall the definition of an A-B slice link (see \cite{FL}, \cite{K1} for a more detailed discussion.)

\begin{definition} \label{A-B slice}
A {\em decomposition} of $D^4$ is a pair of compact
codimension zero submanifolds with boundary $A,B\subset D^4$,
satisfying conditions $(1)-(3)$ below. Denote $$\partial^{+}
A=\partial A\cap
\partial D^4, \; \; \partial^{+} B=\partial B\cap \partial D^4,\; \;
\partial A=\partial^{+} A\cup {\partial}^{-}A, \; \; \partial
B=\partial^{+} B\cup {\partial}^{-}B.$$
(1) $A\cup B=D^4$,\\
(2) $A\cap B=\partial^{-}A=\partial^{-}B,$ \\
(3) $S^3=\partial^{+}A\cup \partial^{+}B$ is the standard genus $1$
Heegaard decomposition of $S^3$.

Given an $n-$component link $L=(l_1,\ldots,l_n)\subset S^3$, let
$D(L)=(l_1,l'_1,\ldots, l_n,l'_n)$ denote the $2n-$component link
obtained by adding an untwisted parallel copy $L'$ to $L$. The link
$L$ is {\em $A-B$ slice} if there exist decompositions $(A_i, B_i),
i=1,\ldots, n$ of $D^4$ and self-homeomorphisms ${\phi}_i, {\psi}_i$
of $D^4$, $i=1,\ldots,n$ such that all sets in the collection
${\phi}_1 A_1, \ldots, {\phi}_n A_n, {\psi}_1 B_1,\ldots, {\psi}_n
B_n$ are disjoint and satisfy the boundary data:
${\phi}_i({\partial}^{+}A_i)$  is a tubular neighborhood of $l_i$
and ${\psi}_i({\partial}^{+}B_i)$ is a tubular neighborhood of
$l'_i$, for each $i$.
\end{definition}

The surgery conjecture is equivalent to the statement that the
Borromean rings, and a certain family of their generalizations, are $A-B$
slice. In \cite{K1} we constructed a decomposition $D^4=A\cup B$ and
disjoint embeddings
$A_i, B_i$ into $D^4$ so that the
attaching curves $\{ {\alpha}_i\}$ of the $A_i$ form the Borromean rings
(or more generally any given link with trivial linking numbers)
and the curves $\{ {\beta}_i \}$ form an untwisted
parallel copy. The validity of one of the conditions necessary for solving the canonical
surgery problems was unknown at the time of that construction, namely the equivariance
(the existence of the homeomorphisms ${\phi}_i, {\psi}_i$), or phrased differently
it was not known whether there exist disjoint re-embeddings of the submanifolds $A, B$ which are
{\em standard}. It follows from theorem \ref{thm} that standard disjoint embeddings for these decompositions
do not exist.
Therefore an open question, important in search for an obstruction to surgery in the context
of the A-B slice problem, is: {\em Given any decomposition $D^4=A\cup B$, is one of the
two embeddings $A\hookrightarrow D^4$, $B\hookrightarrow D^4$ necessarily robust?}

\bigskip

\section{Construction of the submanifolds} \label{decompositions}

\begin{figure}[ht]
\vspace{.5cm}
\includegraphics[width=3.8cm]{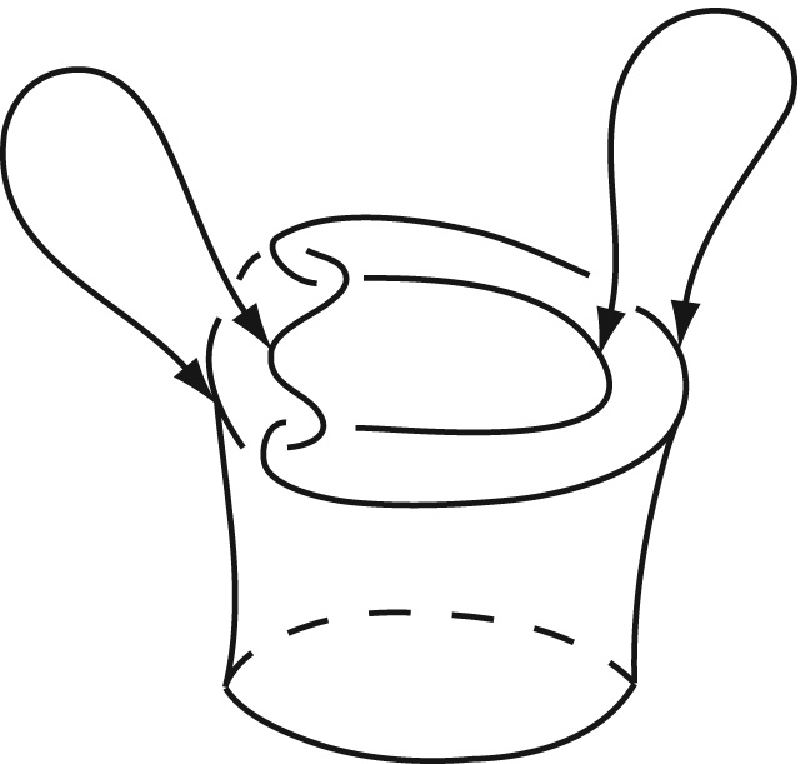} \hspace{3cm} \includegraphics[width=4cm]{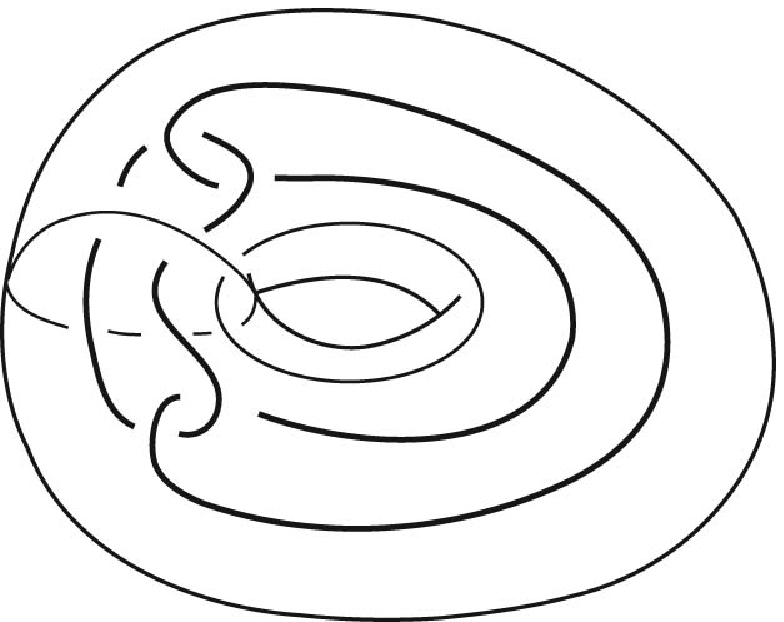}
{\small
    \put(-237,-7){${\beta}$}
    \put(-331,73){$H_1$}
    \put(-205,80){$H_2$}
    \put(-310,23){$B_0$}}
{\scriptsize
    \put(-107,30){$0$}
    \put(-17,22){$0$}
    \put(-46,35){$\beta$}}
\vspace{.45cm} \caption{}
\label{fig:Bing pair}
\end{figure}

This section reviews the construction \cite{K1} of the submanifolds
of $D^4$, which will be used in the proof of theorem \ref{thm} in section \ref{sec:Milnor}.
The construction consists of a series of modifications of the handles
structures, starting with a standard
surface and its complement in the $4-$ball.
Consider the genus one surface
$S$ with a single boundary component ${\alpha}$, and set $A_0=S\times D^2$.
Consider the standard embedding $(S,{\alpha})\subset(D^4,S^3)$
(take an embedding of the surface in $S^3$, push
it into the $4-$ball and take a regular neighborhood.) Then $A_0$ is identified
with a regular neighborhood of $S$ in $D^4$.
The complement $B_0$ of $A_0$ in the $4-$ball is obtained from the collar
on its attaching curve, $S^1\times D^2\times I$, by attaching a pair
of zero-framed $2-$handles to the Bing double of the core of the
solid torus $S^1\times D^2\times\{1\}$, figure \ref{fig:Bing pair}.
(See for example \cite{FL} for a proof of this statement.)

Note that a distinguished pair of curves ${\alpha}_1, {\alpha}_2$,
forming a symplectic basis in the surface $S$, is determined as the
meridians (linking circles) to the cores of the $2-$handles $H_1,
H_2$ of $B_0$ in $D^4$. In other words, ${\alpha}_1$,
${\alpha}_2$ are fibers of the circle normal bundles over the cores
of $H_1, H_2$ in $D^4$.

\begin{figure}[ht]
\vspace{.5cm}
\includegraphics[width=3.8cm]{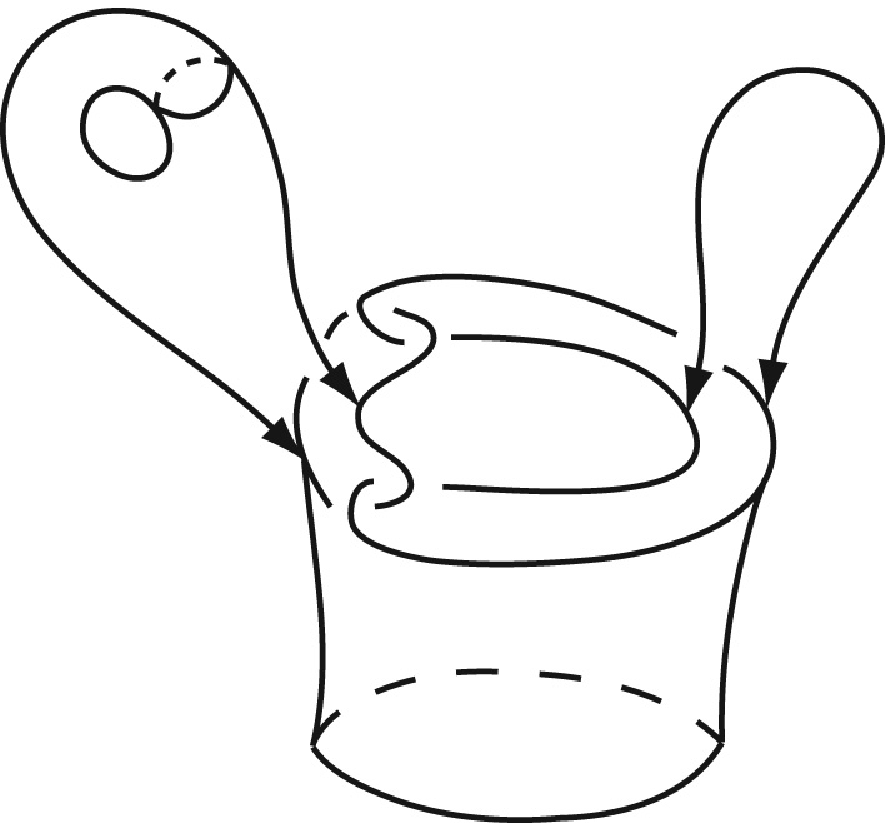} \hspace{2.2cm} \includegraphics[width=6cm]{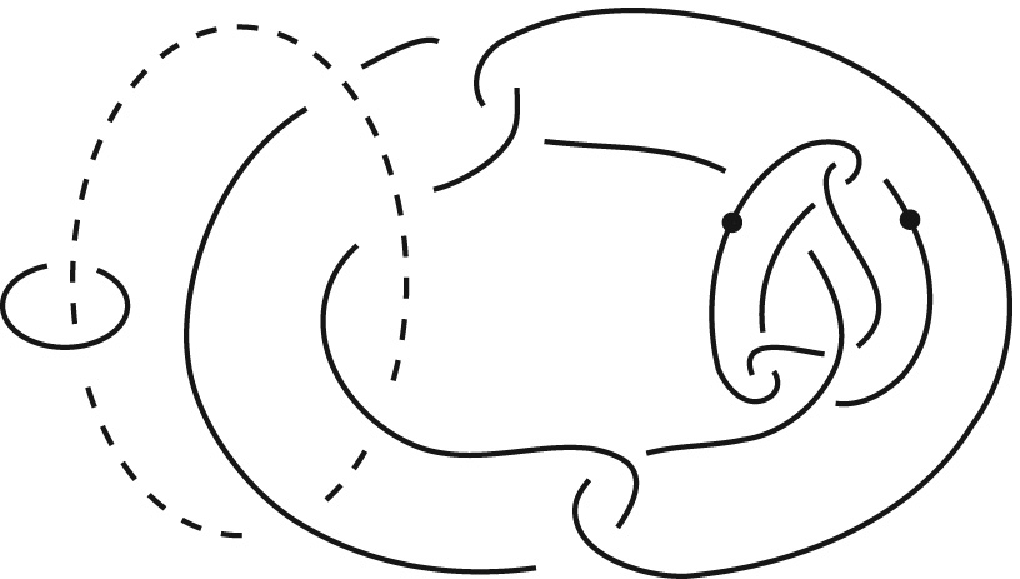}
{\small
    \put(-346,34){$B_1$}
    \put(-330,0){${\beta}$}
    \put(-184,42){${\beta}$}}
{\scriptsize
    \put(-101,-10){$0$}
    \put(-61,-11){$0$}}
    \vspace{.45cm} \caption{}
    \label{fig:B1}
\end{figure}

An important observation \cite{FL} is that this construction may be
iterated: consider the $2-$handle $H_1$ in place of the original
$4-$ball. The pair of curves (${\alpha}_1$, the attaching circle
${\beta}_1$ of $H_1$) form the Hopf link in the boundary of $H_1$.
In $H_1$ consider the standard genus one surface $T$ bounded by
${\beta}_1$. As discussed above, its complement is given by two
zero-framed $2-$handles attached to the Bing double of ${\alpha}_1$.
Assembling this data, consider the new decomposition $D^4=A_1\cup
B_1$ (in this paper we need only the $B-$side of the decomposition, shown in
figure \ref{fig:B1}.) As above, the diagrams are drawn in solid tori
(complements in $S^3$ of unknotted circles drawn dashed in the
figures.) The handlebodies $A_1, B_1$ are examples of {\em model
decompositions} \cite{FL} obtained by iterated applications of the
construction above. It is known that such
model handlebodies are robust, and in particular the Borromean rings
are not weakly $A-B$ slice when restricted to the class of model
decompositions.

\begin{figure}[ht]
\vspace{.7cm}
\includegraphics[width=4.2cm]{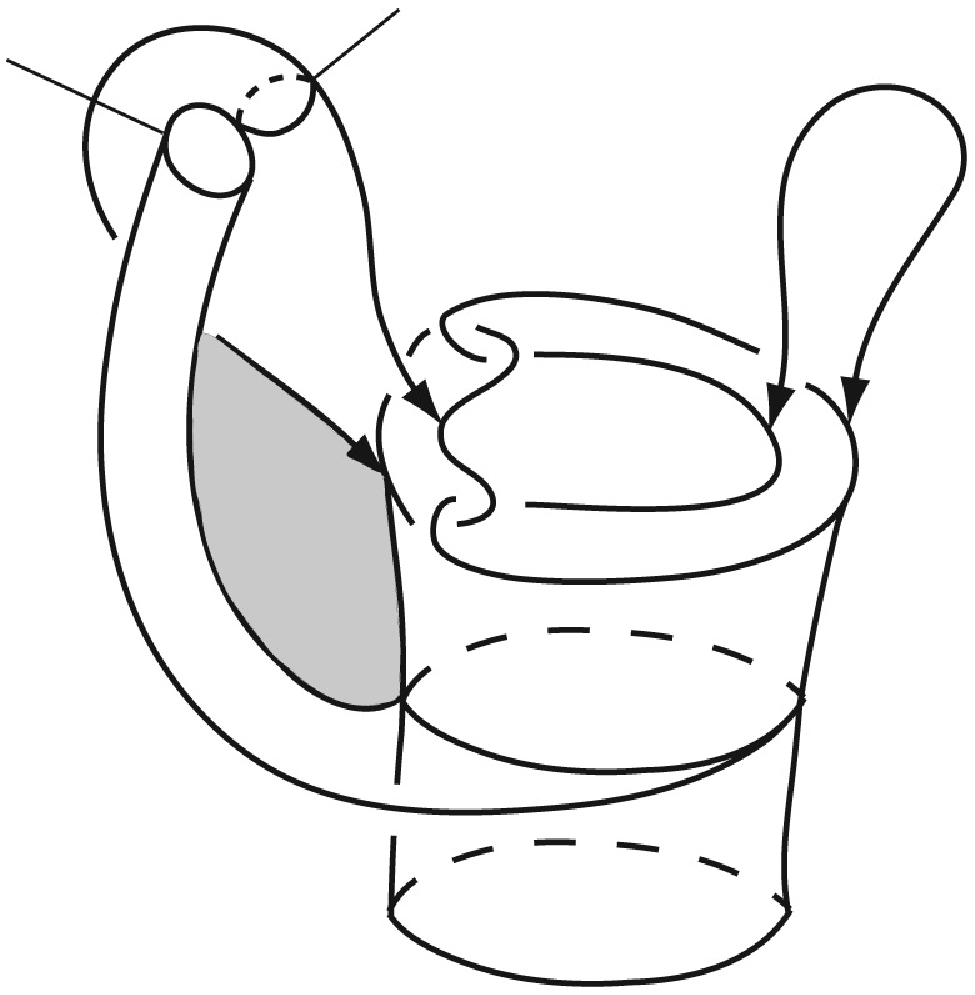} \hspace{2.4cm} \includegraphics[width=6.2cm]{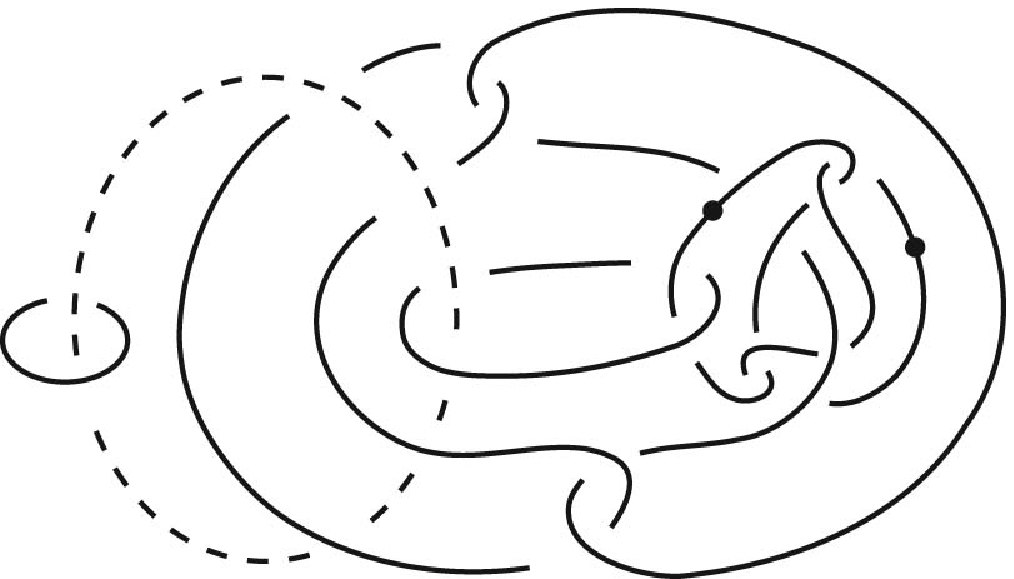}
{\small
    \put(-379,114){$\gamma$}
    \put(-322,120){$\delta$}
    \put(-323,95){$T$}
    \put(-288,100){$D$}
    \put(-366,30){$B$}
    \put(-331,-4){${\beta}$}
    \put(-187,35){$\beta$}}
{\scriptsize
    \put(-98,-10){$0$}
    \put(-58,-11){$0$}
    \put(-80,59){$0$}}
    \vspace{.45cm} \caption{}
\label{fig:B}
\end{figure}

We are now in a position to define the decomposition $D^4=A\cup B$
used in the proof of theorem \ref{thm}.

\begin{definition} \label{definition} Consider
$B=(B_1\, \cup $ zero-framed $2-$handle$)$ attached as shown in the
Kirby diagram in figure \ref{fig:B}.
\end{definition}

Imprecisely (up to homotopy, on the level of spines) $B$ may be
viewed as $B_1\cup 2-$cell attached along (the attaching circle
${\beta}$ of $B_1$, followed by a curve representing a generator of
$H_1$ of the second stage surface of $B_1$). This $2-$cell is
schematically shown in the spine picture of $B$ in the first part of
figure \ref{fig:B} as a cylinder connecting the two curves. The shading
indicates that the new generator of ${\pi}_1$ created by adding the
cylinder is filled in with a disk. The figure showing the
spine is provided only as a motivation for the construction; a
precise description is given by the handle
diagram.

\begin{figure}[ht]
\vspace{.7cm}
\includegraphics[width=6cm]{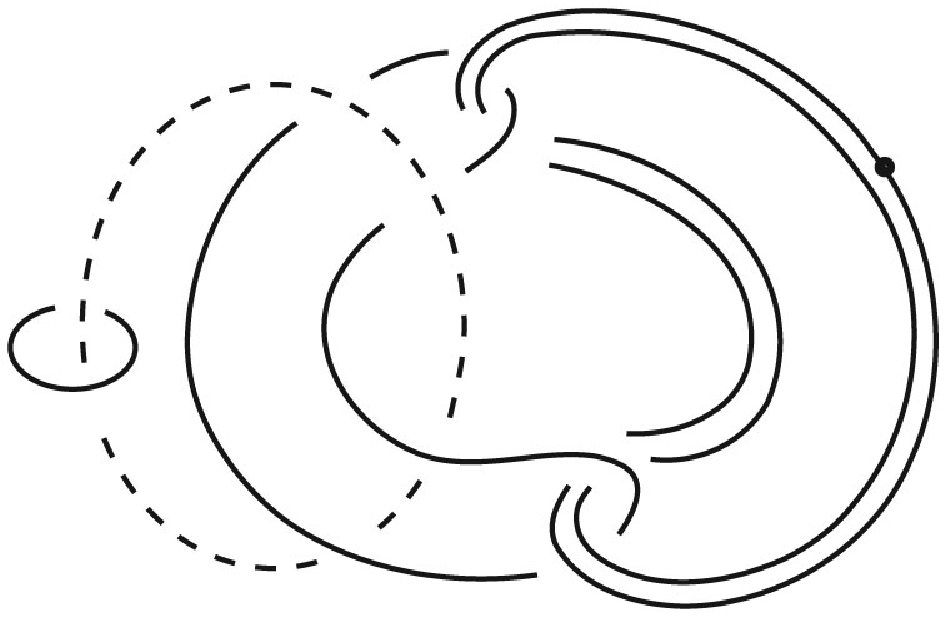}
{\small
    \put(-187,35){$\beta$}}
{\scriptsize
    \put(-90,-4){$0$}
    \put(-48,10){$0$}}
    \vspace{.45cm} \caption{}
\label{fig:Bnew}
\end{figure}

Note that canceling a ($1$-, $2$-)handle pair, one gets the diagram for $B$ shown in figure \ref{fig:Bnew},
this fact will be used in the proof of theorem \ref{thm} in the next section. (Observe that the handle diagram
in figure \ref{fig:Bnew} may also be obtained from the handle diagram of its complement, figure 12 in \cite{K1}.)

\bigskip

\section{The Milnor group and the proof of theorem \ref{thm}} \label{sec:Milnor}

We start this section by summarizing the relevant information about the Milnor group
which will be used in the proof of theorem \ref{thm}. The reader is referred
to the original reference \cite{M} for a more complete introduction to the Milnor
group of links in the $3-$sphere, see also \cite[Section 2]{FT} for a discussion
of the Milnor group in the more general $4-$dimensional context.

\begin{defi}
Given a group $G$, normally generated by elements $g_1,\ldots, g_n$, the
{\em Milnor group} of $G$, relative to the given normal generating set $\{ g_i\}$,
is defined as
\begin{equation} \label{eq:Milnor group}
MG = G \, /\, \langle\! \langle \, [g_i^x, g_i^y] \;\; \forall x,y\in G, i=1,\ldots, n\, \rangle\!\rangle.
\end{equation}
\end{defi}

The Milnor group is a finitely presented nilpotent group of class $\leq n$, where
$n$ is the number of normal generators in the definition above.
An example of interest in this paper is $G={\pi}_1(D^4\smallsetminus {\Sigma})$
where ${\Sigma}$ is a collection of surfaces with boundary, properly and disjointly embedded in
$(D^4, S^3)$. In this case a choice of normal generators is provided by the meridians
$m_i$ to the components ${\Sigma}_i$ of ${\Sigma}$. Here a meridian $m_i$ is an
element of $G$ which is obtained by following a path $\alpha_i$ in $D^4\smallsetminus {\Sigma}$
from the basepoint to the boundary of a regular neighborhood of ${\Sigma}_i$, followed
by a small circle (a fiber of the circle normal bundle) linking ${\Sigma}_i$, then followed
by ${\alpha}_i^{-1}$.



Denote by $F_{g_1,\ldots, g_n}$ the free group generated by the $\{g_i\}$, and consider the Magnus
expansion
\begin{equation}\label{Magnus}
M\co F_{g_1,\ldots, g_n}\longrightarrow {\mathbb Z}[x_1,\ldots, x_n]
\end{equation}
into the ring
of formal power series in non-commuting variables $\{ x_i\}$, defined by $M(g_i)=1+x_i, \,
M(g_i^{-1})=1-x_i+x_i^2\mp\ldots$ We will keep the same notation for the homomorphism
\begin{equation} \label{MagnusMilnor}
M\co MF_{g_1,\ldots, g_n}\longrightarrow R_{x_1,\ldots,x_n},
\end{equation}
induced by the Magnus expansion,
into the quotient of ${\mathbb Z}[x_1,\ldots, x_n]$ by the ideal generated by all monomials
$x_{i_1}\cdots x_{i_k}$ with some index occuring at least twice. It is established in \cite{M} that
the homomorphism (\ref{MagnusMilnor}) is well-defined and injective.

We now turn to the proof of theorem \ref{thm}.
Consider the submanifold $i\co (B,{\beta})\subset (D^4,S^3)$ constructed in definition \ref{definition}.
Part (i) of theorem \ref{thm} follows from \cite[Theorem 1]{K1} which showed that there exist
disjoint embeddings of three copies $(B_i,{\beta}_i)$ such that the link formed by the curves
${\beta}_1, {\beta}_2, {\beta}_3$ in the $3-$sphere is the Borromean rings. It is convenient to introduce the following
definition.

\begin{definition} An embedding $j\co (B,{\beta})\hookrightarrow (D^4,S^3)$ is {\em standard}
if there exists an ambient isotopy between $j$ and the original  embedding $i$ constructed
in definition \ref{definition}.
\end{definition}

Examining the proof of theorem 1 in \cite{K1}, one may check that the embeddings
of the $B_i$, giving rise to the Borromean rings on the boundary, constructed there, are not standard.
We will now show that given disjoint {\em standard} embeddings of several copies $(B_i, {\beta}_i)$ into
the $4-$ball, the link formed by the curves ${\beta}_1,\ldots, {\beta}_n$ in the $3-$sphere is
necessarily homotopically trivial.
We will show that the Borromean rings do not bound disjoint standard embeddings of three copies of $(B,{\beta})$.
The Borromean rings case is the most interesting example from the perspective of the A-B slice problem, the case of
other homotopically essential links is proved analogously.

\begin{figure}[ht]
\vspace{.5cm}
\includegraphics[width=7cm]{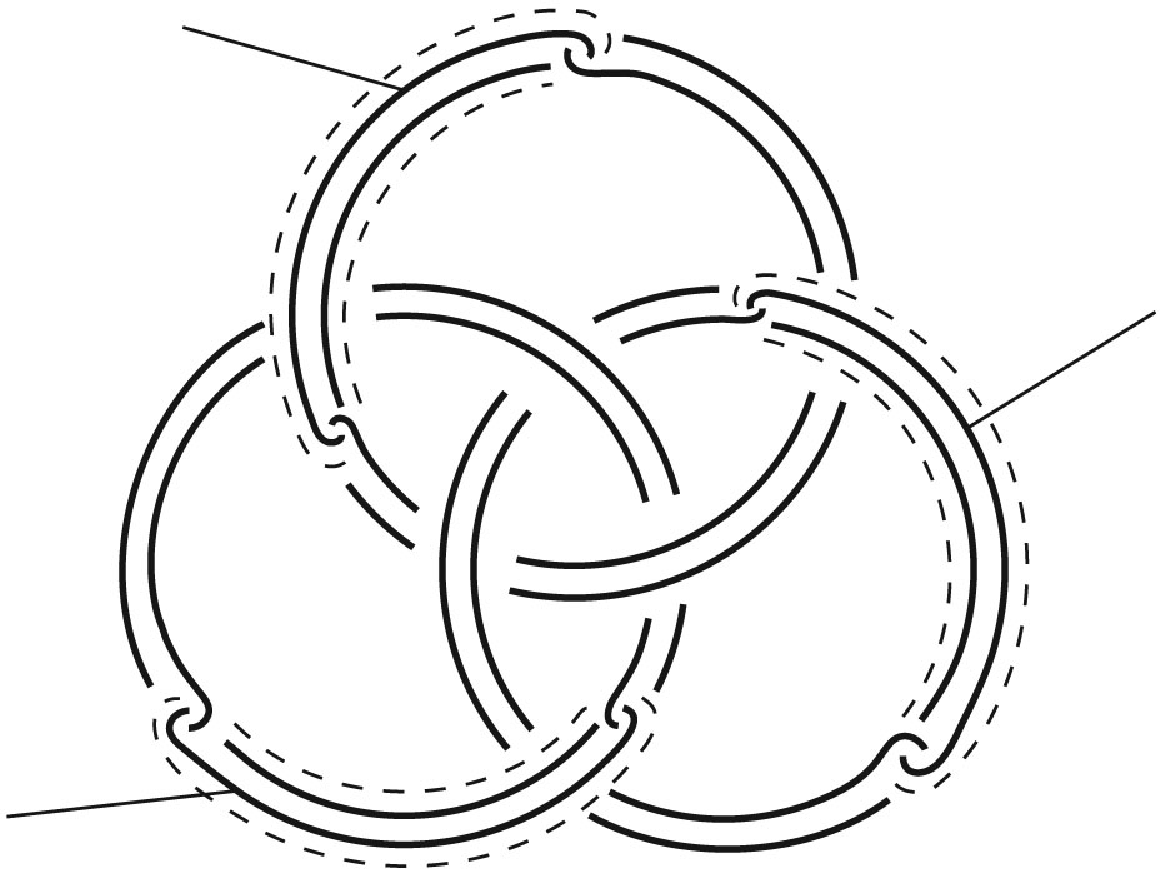}
{\small
    \put(-185,75){$l_1$}
    \put(-180,145){$l_6$}
    \put(-60,130){$l_5$}
    \put(-210,10){$l_2$}
    \put(-77,-9){$l_3$}
    \put(0,95){$l_4$}
    \put(-140,-10){$a$}
    \put(-17,55){$b$}
    \put(-120,152){$c$}}
    \vspace{.45cm} \caption{}
\label{fig:BorDouble}
\end{figure}

Suppose to the contrary that the Borromean rings bound disjoint standard embeddings $B_1,B_2, B_3$.
We will consider the
{\em relative-slice} reformulation of the problem, see \cite{FL} and also \cite{K1} for
a more detailed introduction.
Using the handle diagram in figure
\ref{fig:Bnew}, one observes that then there is a solution to the relative-slice problem shown
in figure \ref{fig:BorDouble}.
This means that the six components $l_1,\ldots,l_6$, drawn solid in the figure,
bound disjoint embedded disks in the
handlebody $D^4\cup_{a,b,c}2$-handles, where the $2$-handles are attached to the $4$-ball with zero framings along
the curves $a,b,c$ drawn dashed in figure \ref{fig:BorDouble}. The fact that the embeddings $B_i\hookrightarrow
D^4$ are standard is reflected in the fact that the slices bounded by the ``solid'' curves of each
$B_i$ do not go over  the $2$-handles (dashed curves) corresponding to the same $B_i$. This means that
the slices for $l_1,l_2$ do not go over $a$, and similarly $l_3, l_4$ do not go over $b$, and $l_5, l_6$ over $c$. (Note that without this restriction, there is a rather straightforward solution to this relative-slice problem.)

Therefore assume the link in figure \ref{fig:BorDouble} is relatively slice, subject to the ``standard'' condition discussed
above. Denote by $D_i$ the slice bounded by $l_i$, $i=1,\ldots, 6$, and let
$D=\cup_{i=2,\ldots, 6} D_i$. Consider $$X:= (D^4\cup_{a,b,c}2{\rm-handles})\smallsetminus D.$$
Denote meridians to the components $l_i$ by $m_i$ and meridians to the curves $a,b,c$ by
$m_a, m_b, m_c$ respectively. The first homology $H_1(X)$ is generated by $m_2,\ldots, m_6$, in fact we
view $\{m_i\}$ as based loops in $X$ normally generating ${\pi}_1(X)$.

Omitting the first component $l_1$, the remaining link $(l_2,\ldots, l_6,a,b,c)$ in figure \ref{fig:BorDouble} is the unlink.
This implies that the second homology $H_2(X)$ is spherical, indeed its generators may be represented by parallel copies
of the cores of the $2$-handles attached to $a,b,c$, capped off by disks in the complement of a neighborhood of the link
$(l_2,\ldots, l_6,a,b,c)$ in the $3$-sphere. Therefore the Milnor group $M{\pi}_1(X)$ with respect to the normal generators $m_i$
is isomorphic to the free Milnor group: $$ M{\pi}_1(X)\cong MF_{m_2,\ldots,m_6}.$$
Indeed, since the Milnor group is nilpotent, it is obtained from the quotient by a term of the lower central series, ${\pi}_1(X)/({\pi}_1(X))^n$, 
by adding the Milnor relations (\ref{eq:Milnor group}). The relations in the nilpotent group ${\pi}_1(X)/({\pi}_1(X))^n$
may be read off from the surfaces representing generators of $H_2(X)$ (see lemma 13 in \cite{K}), in particular the relations
corresponding to spherical classes are trivial. Therefore all relations in $M{\pi}_1(X)$ are the standard relations
(\ref{eq:Milnor group}), or in other words the $M{\pi}_1(X)$ is the free Milnor group.

It follows that the Magnus expansion (\ref{MagnusMilnor}) 
\begin{equation} \label{MagnusMilnor1}
M\co M{\pi}_1(X)\cong MF_{m_2,\ldots,m_6}\longrightarrow R_{x_2,\ldots, x_6}
\end{equation}
is well-defined. Connecting the first component
$l_1$ to the basepoint, consider it as an element of $M{\pi}_1(X)$. According to the assumption the link in figure 
\ref{fig:BorDouble} is relatively slice, therefore $l_1$ bounds a disk in $X$, and in particular it is trivial in 
$M{\pi}_1(X)$. We will find a non-trivial term in the Magnus expansion (\ref{MagnusMilnor1}) 
$M(l_1)\in R_{x_2,\ldots, x_6}$, giving a contradiction with the
relative-slice assumption. 

\begin{figure}[ht]
\vspace{.2cm}
\includegraphics[width=5.5cm]{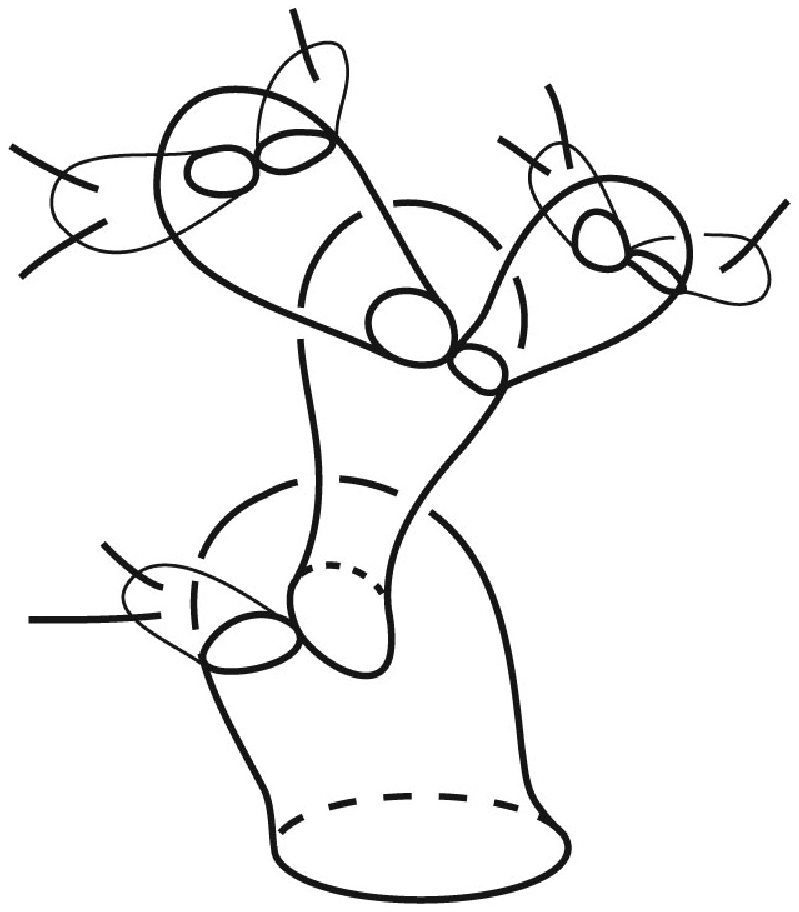}
    \put(-42,5){$l_1$}
    \put(-160,49){$l_2$}
    \put(-144,75){$a$}
    \put(-161,115){$l_4$}
    \put(-162,151){$b$}
    \put(-107,180){$l_3$}
    \put(-69,154){$l_6$}
    \put(-55,165){$c$}
    \put(-7,143){$l_5$}
    \vspace{.2cm} \caption{}
\label{fig:grope}
\end{figure}

Consider meridians $m_i$ to the components $l_i$ in $S^3$, $i=2,\ldots, 6$, and also consider meridians $m_a, m_b, m_c$ to $a,b,c$.
The meridians $m_i$ will also serve as meridians to the slices $D_i$ bounded by $l_i$, $i=2,\ldots,6$ discussed above.
Consider $l_1\in M{\pi}_1(X)$:
\begin{equation} \label{l1}
l_1=[m_a m_2,[[m_3,m_bm_4],[m_5,m_6m_c]]].
\end{equation}
In this expression, $m_a, m_b, m_c$ are elements of $M{\pi}_1(X)$ depending on how the hypothetic slices $D_i$
go over the $2$-handles attached to $a,b,c$.
The expression (\ref{l1}) may be be read off from the capped grope (see figure \ref{fig:grope})
bounded by $l_1$ in the complement of the other components in the $3$-sphere. (Note that the components
$l_2,\ldots,l_6,a,b,c$ intersect only the caps and not the body of the grope.)

Recall a basic commutator identity:  any three elements $f,g,h$ in a group satisfy
\begin{equation} \label{conjugate}
[fg,h]=[f,h]^g\; [g,h].
\end{equation}
Suppose two elements $s,t\in M{\pi}_1(X)$ have Magnus expansions $M(s)=1+{\mathbf x},M(t)=1+{\mathbf y}$,
where ${\mathbf x}, {\mathbf y}$ denote the sum of all monomials of non-zero degree in the expansions of $s,t$.
Then the expansion of the conjugate $tst^{-1}$ is of the form $1+{\mathbf x}+\mathbf{xy}+\mathbf{yx}+\ldots$
In particular, any first non-trivial term in the expansion $M(s)$ also appears in the expansion of any conjugate of $s$, 
$M(s^t)$. The expression for $l_1$ (\ref{l1}) is a $5$-fold commutator, and since the ring $R_{x_2,\ldots,x_6}$
is defined in terms of {\em non-repeating} variables, this implies that any monomial of non-zero degree
in the expansion $M(l_1)$ contains all of the variables $x_2,\ldots x_6$. This means that to read off the Magnus expansion,
any conjugation that comes up while using the identity (\ref{conjugate}) to simplify the expression (\ref{l1}) may be omitted.
These observations imply that the Magnus expansion $M(l_1)$ equals 
\begin{equation} \label{one term}
M([m_2,[[m_3,m_4],[m_5,m_6]]])
\end{equation}
times the Magnus expansion of seven other terms where some (or all) of $m_2,m_4,m_6$ are replaced with $m_a,m_b,m_c$. 
Moreover, recall that the ``standard'' embedding assumption implies that the slice $D_2$ bounded by $l_2$ does not go over
the $2$-handle attached to $a$. Therefore the Magnus expansion of $m_a$ is of the form $M(m_a)=1+\sum_{i=3}^6{\alpha}_i x_i+$higher terms,
for some coefficients ${\alpha}_i$. Since the meridians $m_3, m_5$ are present in each commutator obtained by simplifying 
(\ref{l1}), the only terms in the Magnus expansion of $m_a$ that may contribute to a non-trivial monomial in $M(l_1)$ are
$x_4$ and $x_6$. Similarly, the only possibly non-trivial contributions to $M(l_1)$ of $m_b$ are $x_2, x_6$ and of $m_c$ are 
$x_2, x_4$.

Using the fact that $M([s,t])=1+\mathbf{xy}-\mathbf{yx}\pm\ldots$, where $M(s)=1+{\mathbf x},M(t)=1+{\mathbf y}$,
note that the expansion (\ref{one term}) contains the monomial
$x_2x_3x_4x_6x_5$. We claim that this monomial does not cancel with any other term in the expansion $M(l_1)$.
This claim is proved by a direct inspection: any monomial in the Magnus expansion of a commutator of the form (\ref{one term})
with $m_2$ replaced by $m_a$ has $x_4$ or $x_6$ as either the first or last variable. The only other possibility is 
the expansion of the commutator $[m_2,[[m_3,m_b],[m_5,m_c]]]$ with $M(m_b)$ contributing $x_6$ and $M(m_c)$ contributing
$x_4$. The monomial $x_2x_3x_4x_6x_5$ does not appear in this expansion either. Therefore we found a non-trivial
term in the Magnus expansion $M(l_1)\in M{\pi}_1(X)$, contradicting the relative-slice assumption. 
This contradiction completes the proof.


\begin{thebibliography}{10}

\bibitem{FL} M. Freedman and X.S. Lin, {\em On the $(A,B)$-slice
problem}, Topology Vol. 28 (1989), 91-110.

\bibitem{FQ} M. Freedman and F. Quinn, {\em The topology of
4-manifolds}, Princeton Math. Series 39, Princeton, NJ, 1990.

\bibitem{FT} M. Freedman and P. Teichner,
{\em$4$-Manifold Topology II: Dwyer's filtration and surgery kernels},
Invent. Math. 122 (1995), 531-557.

\bibitem{K} V. Krushkal, {\em Additivity properties of Milnor's $\bar\mu$-invariants},
J. Knot Theory Ramifications 7 (1998), 625-637.

\bibitem{K1} V. Krushkal, {\em A counterexample to the strong version of
Freedman's conjecture}, Ann. of Math. 168 (2008), 675-693.
[arXiv:math/0610865]

\bibitem{KT} V. Krushkal and P. Teichner, {\em Alexander duality,
gropes and link homotopy}, Geom. Topol. 1 (1997), 51-69.
[arXiv:math/9705222]

\bibitem{M} J. Milnor, {\em Link Groups}, Ann. Math 59 (1954), 177-195.



\end{thebibliography}
\end{document}